\documentclass[11pt]{article}
\usepackage{graphicx}
\usepackage{epsfig}
\usepackage{latexsym}
\newsavebox{\toy}
\savebox{\toy}{\framebox[0.65em]{\rule{0cm}{1ex}}}
\newcommand{\QED}{\usebox{\toy}\end{demo}}
\newenvironment{property}%
{\begin{list}{}{\setlength{\rightmargin}{0pt}%
\setlength{\itemsep}{0pt}}}{\end{list}}
\newlength{\templength}
\newcommand{\bp}{\setlength{\templength}{\labelwidth}%
\setlength{\labelwidth}{2em}\begin{property}}
\newcommand{\ep}{\end{property}\setlength{\labelwidth}{\templength}}
\newtheorem{theorem}{Theorem}[subsection]
\newtheorem{lemma}[theorem]{Lemma}
\newtheorem{proposition}[theorem]{Proposition}
\newtheorem{corollary}[theorem]{Corollary}
\newtheorem{assumption}{Assumption}
\newtheorem{definition}{Definition}[subsection]
\newtheorem{remark}{Remark}[subsection]
\newtheorem{exercise}{Exercise}[subsection]

\newcommand{\Thm}[1]{Theorem \ref{Thm.#1}}
\newcommand{\Lem}[1]{Lemma \ref{Lem.#1}}

\newcommand{\Prop}[1]{Proposition \ref{Prop.#1}}

\newcommand{\Theorem}[1]{\begin{theorem}\label{Thm.#1}}
\newcommand{\Lemma}[1]{\begin{lemma}\label{Lem.#1}}
\newcommand{\Proposition}[1]{\begin{proposition}\label{Prop.#1}}
\newcommand{\Corollary}[1]{\begin{corollary}\label{Cor.#1}}
\newcommand{\Assumption}[1]{\begin{assumption}\label{Ass.#1}\rm}
\newcommand{\Definition}[1]{\begin{definition}\label{Def.#1}\rm}
\newcommand{\Remark}[1]{\begin{remark}\label{Rem.#1}\rm }
\newcommand{\Exercise}[1]{\begin{exercise}\label{Exe.#1}\rm }

\newcommand{\bd}{\begin{displaymath}}
\newcommand{\ed}{\end{displaymath}}
\newcommand{\bdn}{\begin{equation}}
\newcommand{\bdnl}{\begin{equation}\label}
\newcommand{\edn}{\end{equation}}
\newcommand{\barray}{\begin{array}}
\newcommand{\earray}{\end{array}}
\newcommand{\bds}{\begin{description}}
\newcommand{\eds}{\end{description}}
\newcommand{\bitemize}{\begin{itemize}}
\newcommand{\eitemize}{\end{itemize}}
\newcommand{\benumerate}{\begin{enumerate}}
\newcommand{\eenumerate}{\end{enumerate}}
\newcommand{\btabbing}{\begin{tabbing}}
\newcommand{\etabbing}{\end{tabbing}}
\newcommand{\bcenter}{\begin{center}}
\newcommand{\ecenter}{\end{center}}
\newcommand{\bflushright}{\begin{flushright}}
\newcommand{\bflushleft}{\begin{flushleft}}
\newcommand{\eflushright}{\end{flushright}}
\newcommand{\eflushleft}{\end{flushleft}}
\newcommand{\bdnn }{\begin{eqnarray*}}
\newcommand{\ednn }{\end{eqnarray*}}
\newcommand{\bdmn}{\begin{eqnarray}}
\newcommand{\edmn}{\end{eqnarray}}

\newcommand{\SSC}[1]{\section{#1}\setcounter{equation}{0}}

\newcounter{biblio}
\newenvironment{references}%
{\begin{list}{[\arabic{biblio}]}{\usecounter{biblio}%
\setlength{\leftmargin}{2.5em}\setlength{\rightmargin}{0pt}%
\setlength{\labelwidth}{2em}\setlength{\itemsep}{0pt}}}{\end{list}}
\newcommand{\References}%
{\vspace{2.8ex plus .3ex minus .3ex}%
\begin{center}{\bf References}\end{center}\begin{references}}

\usepackage{amssymb}
\usepackage{mathrsfs}
\newcommand{\bL}{{\mathbb{L}}}
\newcommand{\N}{{\mathbb{N}}}
\newcommand{\Z}{{\mathbb{Z}}}
\newcommand{\zd}{\Z^d}

\newcommand{\bP}{{\mathbb{bP}}}
\newcommand{\R}{{\mathbb{R}}}

\newcommand{\lra }{\longrightarrow }

\newcommand{\ov}{\overline}
\newcommand{\tl}{\widetilde}

\newcommand{\vvs}{\vspace{2ex}}
\newcommand{\vs}{\vspace{1ex}}

\newcommand{\lan}{\langle \:}
\newcommand{\ran}{\: \rangle}
\newcommand{\lef}{\left}
\newcommand{\rig}{\right}
\newcommand{\ri}{\right}
\newcommand{\st}{\stackrel}

\newcommand{\8}{\infty}

\newcommand{\dps}{\displaystyle}
\newcommand{\sub}{\subset}

\def\P{{\bf P}}

\newcommand{\inflim}{\mathop{\underline{\lim}}}
\newcommand{\suplim}{\mathop{\overline{\lim}}}
\newcommand{\epty}{\emptyset}
\renewcommand{\a}{\alpha}

\newcommand{\del}{\delta}
\newcommand{\D}{\Delta}

\newcommand{\z}{\zeta}
\newcommand{\h}{\eta}

\newcommand{\n}{\nu}

\newcommand{\rh}{\rho}

\newcommand{\s}{\sigma}

\newcommand{\vp}{\varphi}

\newcommand{\W}{\Omega}

\newcommand{\cF }{{\cal F}}

\newcommand{\cP }{{\cal P}}

\newcommand{\cR }{{\cal R}}

\newcommand{\ovn}{\ov{N}}
\newcommand{\ovN}{\ov{N}}

\makeatletter
\def\section{\@startsection{section}{1}{\z@}{-3.5ex plus -1ex minus
 -.2ex}{2.3ex plus .2ex}{\bf}}
\def\subsection{\@startsection{subsection}{2}{\z@}{-3.25ex plus -1ex minus
 -.2ex}{1.5ex plus .2ex}{\bf}}
\makeatother

\def\R{{\cal R}}

\def\P{{P}}
\def\pp{{\cal P}}
\def\z{{\mathbb Z}}

\begin{document}
\parindent=0pt

\bcenter

\large{\bf Localization for Branching Random Walks in Random Environment}

\vvs \normalsize

\noindent Yueyun Hu and Nobuo Yoshida\footnote{Corresponding 
autor, Partially supported by JSPS Grant-in-Aid for Scientific
Research, Kiban (C) 17540112}\\
\medskip {\it  Universit\'e Paris 13 and Kyoto University}

  \ecenter

\begin{abstract}
We consider branching random walks in $d$-dimensional integer
lattice with time-space i.i.d. offspring distributions. This model
is known to exhibit a phase transition: If $d \ge 3$ and the
environment is ``not too random", then, the total population grows
as fast as its expectation with strictly positive probability. If,
on the other hand, $d \le 2$, or the environment is ``random
enough", then the total population grows strictly slower than its
expectation almost surely. We show the equivalence between the
slow population growth and a natural localization property in
terms of ``replica overlap". We also prove a certain stronger
localization property, whenever the total population grows
strictly slower than its expectation almost surely.
\end{abstract}


\normalsize
\SSC{Introduction}
\subsection{Branching random walks in random environment (BRWRE)}
 \label{sec:brw}
We begin by introducing the model. We write $\N=\{0,1,2,...\}$,
$\N^*=\{1,2,...\}$ and $\Z=\{ \pm x \; ; \; x \in \N \}$ in the
sequel. Let $p(\cdot, \; \cdot)$ be a transition probability for
the symmetric simple random walk on $\zd$: \bdnl{1/2d}
p(x,y)=\left\{ \barray{ll}
\frac{1}{2d} & \mbox{if $|x-y|=1$,} \\
0 & \mbox{if $|x-y|\neq 1$,}\earray \rig. \edn where
$|x|=(|x_1|^{2}+..+|x_d|^{2})^{1/2}$ for $x \in \zd$. To each $(t,x)
\in \N \times \zd$, we associate a distribution
$$
q_{t,x}=(q_{t,x}(k))_{k \in \N} \in [0,1]^\N, \; \; \; \sum_{k \in
\N}q_{t,x}(k)=1
$$
on $\N$. Then, the branching random walk (BRW) with offspring
distribution $q=(q_{t,x})_{(t,x) \in \N \times \zd}$ is described
as the following dynamics:

\bitemize \item At time $t=0$, there is one particle at the origin
$x=0$. \item Suppose that there are $N_{t,x}$ particles at each
site $x \in \zd$ at time $t$. At time $t+1$, the $\n$-th particle
at a site $x$ ($\n =1,..,N_{t,x}$) jumps to a site $y=X^\n_{t,x}$
with probability $p(x,y)$ independently of each other.
 At arrival, it dies, leaving $K^\n_{t,x}$ new particles there.
\eitemize

We formulate the above description more precisely.
The following formulation is an analogue of \cite[section
4.2]{Rev94}, where non-random offspring distributions are
considered. See also \cite[section 5]{BGK05} for the random
offspring case.

\vvs \noindent $\bullet$ {\it Spatial motion:} A particle at
time-space location $(t,x)$ is supposed to jump to some other location 
$(t+1,y)$ and is replaced by its children there. Therefore, the
spactial motion should be described by assignning destination
of the each particle at each time-space location $(t,x)$. So, we
are guided to the following definition. We define the measurable
space $(\W_X,\cF_X)$ as the set $(\zd)^{\N \times \zd \times
\N^*}$ with the product $\s$-field, and $\W_X \ni X \mapsto
X^\n_{t,x}$ for each $(t,x,\n) \in \zd \times \N \times \N^*$ as
the projection. We define $P_X \in \cP (\W_X,\cF_X)$ as the
product measure such that \bdnl{P_X} P_X (X^\n_{t,x}=y)=p(x,y)\;
\; \; \mbox{for all $(t,x,\n) \in \N \times \zd \times \N^*$ and
$y \in \zd$.} \edn Here, we interpret $X^\n_{t,x}$ as the position
at time $t+1$ of the children born from the $\n$-th particle at
time-space location $(t,x)$.

\vvs \noindent $\bullet$ {\it Offspring distribution:} We set
$\W_q =\cP (\N)^{\N \times \zd}$, where $\cP (\N)$ denotes the set
of probability measures on $\N$:
$$
\cP (\N)=\{ q=(q (k))_{k \in \N} \in [0,1]^\N \; ; \; \sum_{k \in
\N}q (k)=1\}.
$$
Thus, each $q \in \W_q$ is a function $(t,x) \mapsto
q_{t,x}=(q_{t,x}(k))_{k \in \N}$ from $\N \times \zd$ to $\cP
(\N)$. We interpret $q_{t,x}$ as the offspring distribution for
each particle which occupies the time-space location $(t,x)$. The
set $\cP(\N)$ is equipped with the natural Borel $\s$-field
induced from that of $[0,1]^\N$. We denote by $\cF_q$ the product
$\s$-field on $\W_q$.

We define the measurable space $(\W_K,\cF_K)$ as the set $\N^{\N
\times \zd \times \N^*}$ with the product $\s$-field, and $\W_K
\ni K \mapsto K^\n_{t,x}$ for each $(t,x,\n) \in \N \times \zd
\times \N^*$ as the projection. For each fixed $q \in \W_q$, we
define $P^q_K \in \cP (\W_K,\cF_K)$ as the product measure such
that \bdnl{P^q_K} P^q_K (K^\n_{t,x}=k)=q_{t,x}(k)\; \; \;
\mbox{for all $(x,t,\n) \in \zd \times \N \times \N^*$ and $k \in
\N$.} \edn We interpret $K^\n_{t,x}$ as the number of the children
born from the $\n$-th particle at time-space location $(t,x)$.

We now define the branching random walk in random environment. We
fix a product measure $Q \in \cP (\W_q,\cF_q)$, which describes
the i.i.d. offspring distribution assigned to each time-space
location. Finally, we define $(\W, \cF)$ by
$$
\W=\W_X \times \W_K \times \W_q, \; \; \; \cF=\cF_X \otimes \cF_K
\otimes \cF_q,
$$
and $P^q, P \in \cP (\W, \cF)$ by
$$
P^q=P_X \otimes P^q_K \otimes \del_q, \; \; \; P=\int Q (dq)P^q.
$$
We denote by $N_{t,x}$ the population at time-space location
$(t,x) \in \N \times \zd$, which is defined inductively
 by
 \bdnl{n(x,t)}
N_{0,x}=\del_{0,x},\; \; N_{t,x}=\sum_{y \in \zd}\sum_{\n
=1}^{N_{t-1,y}} \del_x (X^\n_{t-1,y})K^\n_{t-1,y}, \; \; t \ge 1.
 \edn
We consider the filtration: \bdnl{cF_t} \cF_0=\{ \epty, \W\},\; \;
\cF_t =\s (X^\cdot_{s, \cdot}, K^\cdot_{s, \cdot}, q_{s, \cdot}\;
; \; s \le t-1 )\; \; \; t \ge 1, \edn which the process $t
\mapsto (N_{t,x})_{x \in \zd}$ is adapted to. The total population
at time $t$ is then given by \bdnl{n_t} N_t =\sum_{x \in
\zd}N_{t,x}=\sum_{y \in \zd} \sum_{\n =1}^{N_{t-1,y}}K^\n_{t-1,y}.
\edn We remark that the total population is exactly the classical
Galton-Watson process if $q_{t,x} \equiv q$, where $q \in \cP
(\N)$ is non-random. On the other hand, if $\zd$ is replaced a
singleton, then $N_t$ is the polulation of the Smith-Wilkinson
model \cite{SW69}.

\hspace{3mm} For $p >0$, we write \bdmn m^{(p)}&=&
Q[m^{(p)}_{t,x}] \; \; \mbox{with} \; \;
m^{(p)}_{t,x}=\sum_{k \in \N}k^pq_{t,x}(k), \label{m_p}\\
m &=& m^{(1)}. \label{m_1} \edmn Note that for $p \ge 1$,
$$
m^p \le Q[m_{t,x}^p] \le m^{(p)}
$$
by H\"older's inequality.
We set \bdnl{ovn_t} \ovn_{t,x}=N_{t,x}/m^t\; \; \mbox{and}\; \;
\ovn_t=N_t/m^t. \edn $\ovn_t=N_t/m^t$ is a martingale, and
therefore the following limit always exists: 
\bdnl{ovn_8}
\ovn_\8=\lim_t \ovn_t, \; \; \mbox{$P$-a.s.} \edn 
We denote the
density of the population by: \bdnl{rh}
\rh_{t,x}=\frac{N_{t,x}}{N_t} =\frac{\ov{N}_{t,x}}{\ov{N}_t}, \;
\; t \in \N, x \in \zd \edn Interesting objects related to the
density would be \bdnl{rh^*} \rh^*_t=\max_{x \in \zd}\rh_{t,x}, \;
\; \mbox{and}\; \; \cR_t=\sum_{x \in \zd}\rh_{t,x}^{2}. \edn
$\rh^*_t$ is the density at the most populated site, while $\cR_t$
is the probability that a given pair of particles at time $t$ are
at the same site. We call $\cR_t$ the {\it replica overlap}, in
analogy with the spin glass theory. Clearly, $(\rh^*_t)^{2} \le \cR_t
\le \rh^*_t$. These quantities convey information on
localization/delocalization of the particles. Roughly speaking,
large values of $\rh^*_t$ or $\cR_t$ indicates that the most of
the particles are concentrated on small numbers of ``favorite
sites" ({\it localization}), whereas small values of them implies
that the particles are spread out over large number of sites ({\it
delocalization}).

\subsection{The phase transition in terms of the population growth}
\label{sec.growth}
Due to the random environment, the population $N_t$ has much more
fluctuation as compared with the non-random environment case,
e.g.,\cite[section 4.2]{Rev94}. This fluctuation results from
``disastrous locations" in time-space, where the offspring
distribution $q_{t,x}(k)$ happens to assign extremely high
probability to small $k$'s. Thanks to the random walk, on the
other hand, some of the particles are lucky enough to avoid those
disastrous locations. Therefore, the spatial motion component of
the model has the effect to moderate the fluctuation, while the
random environment intensifies it. These competing factors in the
model give rise to a phase transition as we discuss below.

We first look at the case where the randomness of the offspring
distribution is well moderated by that of the random walk.

Let $(S_t)$ be two a simple
symmetric random walks on $\Z^d$, starting from $0$.  We
denote  by $\pi_d$ the  probability  of the event 
$\cup_{t \ge1}\{S_t=0 \}$. 
As is well known $\pi_d <1$ if and only if $d \ge 3$.
\Proposition{WD}
\bds \item[(a)] There exists a constant $\a^* >\frac{1}{\pi_d}$
such that, if \bdnl{WD1} m>1,\; \;m^{(2)}<\8,\; \; d \ge 3, \; \;
\mbox{and}\; \; \a \st{\rm def.}{=}\frac{Q[m_{t,x}^{2}]}{m^{2}} <\a^*,
\edn then, $P (\ovN_\8 >0)>0$. \item[(b)] If one assumes the
stronger assumption \bdnl{WD2} m>1,\; \;m^{(2)}<\8,\; \; d \ge 3,
\; \; \mbox{and}\; \; \a  <\frac{1}{\pi_d}, \edn then
$$
\cR_T =O (T^{-d/2})\; \;
\; \mbox{in $P \lef( \cdot | \ovN_\8 >0 \rig)$-probability,}
$$
i.e., the laws $P \lef(T^{d/2}\cR_T \in \cdot | \ovN_\8 >0 \rig)$, $T
\ge 1$ are tight. \eds
\end{proposition}
Conditions (\ref{WD1}) and (\ref{WD2}) control the
randomness of the environment in terms of the random walk.
\Prop{WD}(a) says that, under (\ref{WD1}), the total population
grows as fast as its expectation with strictly positive
probability. This was obtained in \cite[Theorem 4]{BGK05}.
\Prop{WD}(b) is a quantative statement for delocalization under
(\ref{WD2}) in terms of the replica overlap \cite[Proposition
1.2.3]{Yos07}.

Next, we turn to the case where the randomness of the environment
dominates:
\Proposition{SD}
Suppose one of the following conditions: \bds \item[(a1)] $d=1, \;
Q(m_{t,x} = m) \neq 1.$ \item[(a2)] $d=2, \; Q(m_{t,x} = m) \neq
1.$ \item[(a3)] ${\dps d \ge3, \;
Q\lef[\frac{m_{t,x}}{m}\ln\frac{m_{t,x}}{m}\ri]>\ln (2d)}$. \eds
Then, $P(\ovn_\8 =0)=1$. Moreover, in cases (a1) and (a3), there
exists a non-random number $c>0$ such that \bdnl{N_t<m^t}
\suplim_{t}\frac{\ln \ov{N}_t }{t}< -c, \; \; \mbox{a.s.} \edn
\end{proposition}
\Prop{SD} says that the total population grows strictly slower
than its expectation almost surely, in low dimensions or in
``random enough" environment. The result is in contrast with the
non-random environment case, where $P(\ovn_\8 =0)=1$ only for
offspring distributions with very heavy tail, more precisely, if
and only if $P [K_{t,x}^\n \ln K_{t,x}^\n]=\8$ \cite[page 24,
Theorem 1]{AtNe72}. Here, we can have $P(\ovn_\8 =0)=1$ even when
$K_{t,x}^\n$ is bounded. Also, (\ref{N_t<m^t}) is in sharp
contrast with the non-random environment case, where it is well
known --see e.g., \cite[page 30, Theorem 3]{AtNe72} --that
$$
\{ N_\8 >0\}\st{\mbox{\scriptsize a.s.}}{=} \{ \lim_{t}\frac{\ln
\ov{N}_t }{t}=0\}\; \; \; \mbox{whenever $m >1$.}
$$
\Prop{SD} was obtained in \cite[Theorem 4]{BGK05} without
(\ref{N_t<m^t}), and in \cite[Corollary 3.3.2]{Yos07} with
(\ref{N_t<m^t}).
\subsection{The results: the localization/delocalization transition}

 In this paper, we aim at the localization problem for the
branching random walk in random environment.
We shall prove that for $d=1, 2$ and for ``random enough
environment" in $d\ge 3$, almost surely, there exists a sequence
of time $t$'s such that both the maximal density $\rho^*_t$ and
overlap $\R_t$ are  bigger than some positive constant.

\vvs We first characterize the event $\{ \ov{N}_\8 =0\}$ in terms
of the replica overlap. Thanks to this characterization, we can
rigorously identify the phase transition in terms of population
growth as discussed in section \ref{sec.growth} with the
localization/delocalization transition in terms of the replica
overlap.
\Theorem{loc}
Suppose that \bdnl{1<K<K}
 m^{(3)}<\8,
\; \; Q(m_{t,x} = m) \neq 1, \; \; Q(q_{t,x}(0)=0)=1. \edn
Then, \bdnl{loc1} \{ \ov{N}_\8 =0 \}\st{\scriptsize
\mbox{a.s.}}{=} \{ \sum_{s=0}^{\8}\cR_s =\8 \}, \edn where
$(\cR_t)_{t \ge 0}$ is defined by (\ref{rh^*}). Moreover, there
exist constants $c_1,c_2 \in (0,\8)$ such that, \bdnl{loc2} \{
\ov{N}_\8 =0\}
 \st{\mbox{\scriptsize a.s.}}{\sub}
\{ \; -c_1\ln \ov{N}_t \le \sum_{s=0}^{t-1}\cR_s \le -c_2\ln
\ov{N}_t\; \; \; \mbox{for large enough $t$'s.}\}. \edn
\end{theorem}
We will prove \Thm{loc} in section \ref{sec.loc}.

\vvs As we referred to before, the large values of the replica
overlap, or the maximal density, indicates the localization of the
particles to a small number of sites. We have the following lower
bound for the replica overlap and the maximal density:
\Theorem{sloc}
Suppose (\ref{1<K<K}) and that $P(\ov{N}_\8=0) =1$. Then, there
exists a non-random number $c \in (0,1)$ such that \bdnl{sloc}
\suplim_{t \nearrow \8}\rh^*_t \ge \suplim_{t \nearrow \8}\cR_t
\ge c, \; \; \mbox{a.s.,} \edn where $(\rh^*_t)_{t \ge 0}$ and
$(\cR_t)_{t \ge 0}$ are defined by (\ref{rh^*}). In particular,
(\ref{sloc}) holds true if we assume any one of (a1) -- (a3) in
\Prop{SD}.
\end{theorem}
(\ref{sloc}) says that the replica overlap persists,
 in contrast with \Prop{WD}(b),
where the replica overlap $\cR_T$ decays like $O(T^{-d/2})$. The
proof of \Thm{sloc} will be  presented in section \ref{sec.sloc}.
Some more remarks on \Thm{sloc} are in order:

\vs \noindent {\bf 1)}
In cases (a1) and (a3) in \Prop{SD}, (\ref{sloc}) follows easily
from (\ref{N_t<m^t}) and (\ref{loc2}). However, the proof we
present does not rely on (\ref{N_t<m^t}), so that we can cover two
dimensional case (a2) as well.

\vs \noindent {\bf 2)}
We prove (\ref{sloc}) by way of the following stronger estimate:
\bdnl{sloc3/2} \inflim_{t \nearrow \8} {\sum_{s=0}^t\cR_s^{3/2}
\over \sum_{s=0}^t\cR_s} \ge c, \; \; a.s. \edn for some
non-random number $c>0$. This in particular implies the following
quantative lower bound on the number of times, at which the
replica overlap is larger than a certain positive number:
$$
\inflim_{t \nearrow \8} {\sum_{s=0}^t1_{\{ \cR_s \ge \epsilon \}}
\over \sum_{s=0}^t\cR_s}
 \ge \epsilon, \; \; a.s.
$$
for small enough $\epsilon >0$.

\vs \noindent {\bf 3)}
For both \Thm{loc} and \Thm{sloc}, similar results are known for
the directed polymers in random environment (DPRE)
\cite{CH02,CSY03,CSY04}. In fact, we have used ideas and
techniques from the DPRE case. However, the results for DPRE do
not seem to directly imply our results.
\SSC{Proof of \Thm{loc}} \label{sec.loc}
\subsection{Lemmas}
For sequences $(a_t )_{t \in \N}$ and $(b_t )_{t \in \N}$ (random
or non-random), we write $a_t \preceq b_t$ if there exists
non-random constant $c \in (0,\8)$ such that $a_t \le c b_t$ for
all $t \in \N$. We write $a_t \asymp b_t$ if $a_t \preceq b_t$ and
$b_t \preceq a_t$.
\Lemma{<N>}
\bds \item[(a)] If $m^{(2)}<\8$ and $ Q(m_{t,x} = m) \neq 1 $,
then, ${\dps P[(N_t-mN_{t-1})^{2}|\cF_{t-1} ] \asymp \sum_{x \in
\zd}N_{t-1,x}^{2}}$. \item[(b)] If $m^{(3)}<\8$, then  ${\dps \lef|
P[(N_t-mN_{t-1})^{3}|\cF_{t-1} ] \ri| \preceq \sum_{x \in
\zd}N_{t-1,x}^{3}}$. \eds
\end{lemma}
Proof: (a): Since
$$
N_t-mN_{t-1}=\sum_{x}\sum_{\n =1}^{N_{t-1,x}}(K^\n_{t-1,x}-m),
$$
we have $(N_t-mN_{t-1})^{2}=\sum_{x_1,x_2}F_{x_1,x_2}$, where
$$
F_{x_1,x_2} =\sum_{\nu_1 =1}^{N_{t-1,x_1}}\sum_{\nu_2
=1}^{N_{t-1,x_2}} (K^{\nu_1}_{t-1,x_1}-m)(K^{\nu_2}_{t-1,x_2}-m).
$$
If $x_1 \neq x_2$, then $K^{\nu_1}_{t-1,x_1}$ and
$K^{\nu_2}_{t-1,x_2}$ are mean $m$ independent r.v.'s under $P(
\cdot | \cF_{t-1})$, and hence
$$
P[F_{x_1,x_2}|\cF_{t-1} ]=0.
$$
We may therefore focus on the expectation of $F_{x_1,x_2}$ with
$x_1=x_2=x$. In this case, $\{K^\nu_{t-1,x}\}_{\nu=1}^{N_{t-1,x}}$
are independent under $P( \cdot | \tl{\cF}_{t-1})$, where
$$
\tl{\cF}_{t-1}=\s (\cF_{t-1}, (q_{t-1,x})_{x \in \zd}).
$$
Thus,
$$
P[F_{x,x}|\tl{\cF}_{t-1} ] =N_{t-1,x}(N_{t-1,x}-1)(m_{t-1,x}-m)^{2}
+N_{t-1,x}P^q[(K^\nu_{t-1,x}-m)^{2}].
$$
The first and second terms on the right-hand-side come
respectively from off-diagonal and diagonal terms in $F_{x,x}$. We
now set $\a\st{\rm def.}{=}Q[m_{t,x}^{2}]/m^{2}$. Then, $\a>1$ (since
$Q(m_{t,x} = m) \neq 1$) and \bdnn P[F_{x,x}|\cF_{t-1} ]
& = & (\a-1)m^2N_{t-1,x}(N_{t-1,x}-1)+(m^{(2)}-m^{2})N_{t-1,x}\\
& = & (\a-1)m^2N_{t-1,x}^{2}+(m^{(2)}-\a m^{2})N_{t-1,x}. \ednn
Therefore,
$$
P[(N_t-mN_{t-1})^{2}|\cF_{t-1}]= (\a -1)m^{2}\sum_{x}N_{t-1,x}^{2}
+(m^{(2)}-\a m^{2})N_{t-1},
$$
which implies the desired bound. \\
(b): We have $(N_t-mN_{t-1})^{3}=\sum_{x_1,x_2,x_3}F_{x_1,x_2,x_3}$,
where
$$
F_{x_1,x_2,x_3} =\sum_{\n_1 =1}^{N_{t-1,x_1}}\sum_{\nu_2
=1}^{N_{t-1,x_2}} \sum_{\n_3 =1}^{N_{t-1,x_3}}
(K^{\nu_1}_{t-1,x_1}-m)(K^{\nu_2}_{t-1,x_2}-m)(K^{\nu_3}_{t-1,x_3}-m).
$$
If, for example, $x_1 \not\in \{x_2,x_3 \}$, then
$K^{\nu_1}_{t-1,x_1}$ is independent of $\{ K^{\nu_2}_{t-1,x_2},
K^{\nu_3}_{t-1,x_3}\}$ under $P( \cdot | \cF_{t-1})$, and hence
$P[F_{x_1,x_2,x_3}|\cF_{t-1}]=0$. This implies that
$$
P[(N_t-mN_{t-1})^{3}|\cF_{t-1}]=\sum_{x}P[F_{x,x,x}|\cF_{t-1}].
$$
We have on the other hand that, \bdnn P[F_{x,x,x}|\tl{\cF}_{t-1}]
& = &
N_{t-1,x}P^q[(K^{\nu}_{t-1,x}-m)^{3}] \\
& & +3N_{t-1,x}(N_{t-1,x}-1)
P^q[(K^{\nu}_{t-1,x}-m)^{2}]P^q[K^{\nu}_{t-1,x}-m] \\
& & +N_{t-1,x}(N_{t-1,x}-1)(N_{t-1,x}-2) P^q[K^{\nu}_{t-1,x}-m]^{3}.
\ednn and therefore that
$$
\lef| P[F_{x,x,x}|\tl{\cF}_{t-1}] \ri| \le
N_{t-1,x}^{3}P^q[|K^{\nu}_{t-1,x}-m|^3].
$$
Putting things together, we obtain
$$
\lef| P[(N_t-mN_{t-1})^{3}|\cF_{t-1}] \ri|  \le c\sum_{x}N_{t-1,x}^{3},
\; \; \; \mbox{with $c=Q[|K^{\nu}_{t-1,x}-m|^3]$}.
$$
\hfill $\Box$

\vvs
 Let us now recall Doob's decomposition in our settings.
An $(\cF_t)$-adapted process $X=(X_t )_{t \ge 0} \sub \bL1(P)$ can
be decomposed in a unique way as
$$
X_t=M_t (X)+A_t (X),\; \; \; t \geq 1,
$$
where $M (X)$ is an $(\cF_t)$-martingale and
$$
A_0 =0, \; \; \D A_t =P[\D X_t |\cF_{t-1}], \; \; \; t \ge 1.
$$
Here, and in what follows, we write $\D a_t=a_t-a_{t-1}$ ($t \geq
1$) for a sequence $(a_t )_{t \in \N}$ (random or non-random).
$M_t (X)$ and $A_t (X)$ are called respectively, the martingale
part and compensator of the process $X$. If $X$ is a square
integrable martingale, then the compensator $A_t (X^{2})$ of the
process $X^{2}=( X_t^{2} )_{t \ge 0} \sub \bL1(P)$ is denoted by $\lan X
\ran_t$ and is given by the following formula:
$$
\D \lan X \ran_t =  P[(\D X_t)^{2}|\cF_{t-1}].
$$
Now, we turn to the Doob's decomposition of $X_t=-\ln \ov{N}_t$,
whose martingale part and the compensator will be henceforth
denoted $M_t$ and $A_t$ respectively; \bdnl{dec} -\ln \ov{N}_t
=M_t+A_t, \; \; \; \D A_t=-P[\D \ln \ov{N}_t | \cF_{t-1}] \edn
\Lemma{<M>}
Suppose (\ref{1<K<K}). Then, $\D \lan M \ran _t \preceq \cR_{t-1}
\asymp \D A_t$.
\end{lemma}
Proof: We set $U_t =\frac{\D \ov{N}_t}{\ov{N}_{t-1}}$ to simplify
the notation. We first note the following: \bds \item[(1)] ${\dps
U_t \ge  \frac{1}{m}-1 >-1}$. \item[(2)] ${\dps |\D \ln \ov{N}_t|
\le m\lef| U_t \ri|}$. \item[(3)] ${\dps P \lef[ U_t^{2} | \cF_{t-1}
\ri] \asymp P \lef[ \vp \lef( U_t \ri) | \cF_{t-1} \ri] \asymp
\cR_{t-1}}$, where $\vp (x)=x-\ln (1+x)$. \eds In fact, $N_{t-1}
\le N_t$ by (\ref{1<K<K}), and hence $(1/m)\ov{N}_{t-1} \le
\ov{N}_t$. These imply (1). (2) follows directly from  (1) since
$$
|\ln x -\ln y| \le \frac{m|x-y|}{y} \; \; \; \mbox{if $x,y >0$ and
$x/y \ge 1/m$}.
$$
As for (3), we have by \Lem{<N>}(a) that
$$
P \lef[ U_t^{2} | \cF_{t-1} \ri]
 =  \frac{P \lef[ | \D \ov{N}_t|^2 | \cF_{t-1} \ri]}{\ov{N}_{t-1}^{2}}
\asymp \cR_{t-1}.
$$
We now note that there exists $c \in (0,\8)$, which depends only
on $m$ such that
$$
\frac{x^{2}}{4(2+x)} \le \vp (x) \le cx^{2} \; \; \; \mbox{for all $x
\ge \frac{1}{m}-1$.}
$$
This, together with (1) implies that
$$
P \lef[ \vp \lef( U_t \ri) | \cF_{t-1} \ri] \le cP \lef[ U_t^{2} |
\cF_{t-1} \ri]\asymp \cR_{t-1}.
$$
On the other hand, we have by \Lem{<N>}(b) that
$$
\lef|  P \lef[ U_t^{3} | \cF_{t-1} \ri] \ri| =\frac{1}{N_{t-1}^{3}}\lef|
P[(N_t-mN_{t-1})^{3}|\cF_{t-1} ] \ri| \preceq
\frac{1}{N_{t-1}^{3}}\sum_{x \in \zd}N_{t-1,x}^{3} \le \cR_{t-1}.
$$
Thus, \bdnn \cR_{t-1} & \asymp  & P \lef[ U_t^{2} | \cF_{t-1} \ri]
  =  P\lef[ \frac{U_t}{\sqrt{2+U_t}}U_t\sqrt{2+U_t} | \cF_{t-1}\rig] \\
& \le  & P\lef[ \frac{U_t^{2}}{2+U_t} | \cF_{t-1}\rig]^{1/2} P \lef[
2U_t^{2}+U_t^{3} | \cF_{t-1}\rig]^{1/2}
 \preceq  P \lef[ \vp (U_t) | \cF_{t-1}\rig]^{1/2}\cR_{t-1}^{1/2},
\ednn and hence $\cR_{t-1} \preceq P \lef[ \vp (U_t) |
\cF_{t-1}\rig]$.

\hspace{3mm} The rest of the proof is easy. We have by (3) that
$$
\D A_t  =  -P[\D \ln \ov{N}_t | \cF_{t-1}] =-P \lef[ \ln \lef(
1+U_t\ri) | \cF_{t-1} \ri]
 =  P \lef[ \vp \lef( U_t \ri) | \cF_{t-1} \ri] \asymp \cR_{t-1}.
$$
Similarly, by (2) and \Lem{<N>},
$$
P[(\D \ln \ov{N}_t)^{2} | \cF_{t-1}] \preceq P \lef[ U_t^{2} | \cF_{t-1}
\ri] \asymp \cR_{t-1}.
$$
This, together with $\D A_t \asymp \cR_{t-1}$ implies that
$$
\D \lan M \ran _t=P[(\D M_t)^{2} | \cF_{t-1}] \le 2P[(\D \ln
\ov{N}_t)^{2} | \cF_{t-1}] +2(\D A_t)^{2} \preceq \cR_{t-1}.
$$
\hfill $\Box$

\subsection{Proof of \Thm{loc}} \label{pf_loc}
The proof is based on the decomposition (\ref{dec}). It is enough
to prove the following: \bds \item[(1)] $\{ \ov{N}_\8 =0\}
 \st{\mbox{\scriptsize a.s.}}{\sub}
\{ A_\8 =\8 \} = \lef\{ \sum_{s=0}^\8\cR_s =\8 \rig\}$. \item[(2)]
$\lef\{ \sum_{s=0}^\8\cR_s =\8 \rig\}  \st{\mbox{\scriptsize
a.s.}}{\sub} \{\; -c_1\ln \ov{N}_t \le \sum_{s=0}^{t-1}\cR_s \le
-c_2\ln \ov{N}_t\; \; \; \mbox{for large enough $t$'s.} \}$. \eds
To prove these, we recall the following general facts on square
integrable martingales--see for example \cite[page 252, (4.9) and
page 253, (4.10)]{Dur95}: \bds \item[(3)] $\{ \lan M \ran_\8 <\8
\}
 \st{\mbox{\scriptsize a.s.}}{\sub}
\{ \mbox{${\dps \lim_t M_t}$ converges.} \}.$ \item[(4)] ${\dps \{
\lan M \ran_\8 =\8 \}
 \st{\mbox{\scriptsize a.s.}}{\sub}
\{  \lim_t \frac{M_t}{\lan M \ran_t}=0 \}.}$ \eds
By (3) and \Lem{<M>}, we get (1) as follows: \bdnn \lef\{
\sum_{s=0}^{t-1}\cR_s <\8 \rig\} & = & \{ A_\8 <\8 \}
 = \{ A_\8 <\8, \; \lan M \ran_\8 <\8 \} \\
& \st{\mbox{\scriptsize a.s.}}{\sub} & \{ A_\8 <\8, \;
\mbox{${\dps \lim_t M_t}$ converges.} \}
 \sub  \{ \ov{N}_\8 >0\}.
\ednn We now turn to (2). Since $\{ A_\8 =\8 \}= \lef\{
\sum_{s=0}^\8\cR_s =\8 \rig\}$ and
$$
-\frac{\ln \ov{N}_t}{\sum_{s=0}^{t-1}\cR_s} \asymp -\frac{\ln
\ov{N}_t}{A_t}=\frac{M_t}{A_t}+1,
$$
by \Lem{<M>}, (2) is a consequence of: \bds \item[(5)] $\{ A_\8
=\8 \}  \st{\mbox{\scriptsize a.s.}}{\sub}
 \lef\{ \frac{M_t}{A_t} \lra 0 \rig\}$.
\eds Let us  suppose that $A_\8 =\8$. If $\lan M \ran_\8 <\8 $,
then again by (3), ${\dps\lim_t M_t}$ converges and therefore (5)
holds. If, on the contrary, $\lan M \ran_\8 =\8 $, then by (4) and
\Lem{<M>},
$$
\frac{M_t}{A_t}= \frac{M_t}{\lan M \ran_t} \frac{\lan M
\ran_t}{A_t} \lra 0 \; \; \; \mbox{a.s.}
$$
Thus, (5) is true in this case as well.
\hfill $\Box$
\SSC{Proof of \Thm{sloc}} \label{sec.sloc}
We shall prove \Thm{sloc} in the same spirit of that of
\cite{CH02}. In the following subsection, we give  some
preliminary estimates and the final proof is given in the last
subsection.

\subsection{Lemmas}
A technical result at first:
\Lemma{CH}
Let $\h_i$, $1 \le i \le n$ ($n \ge 2$) be positive independent
random variables on a probability space with the probability
measure $\bP$, such that $\bP[\h_i^{3} ] <\8$ for $i=1,..,n$.Then,
\bdmn \bP\left[\frac{\h_1\h_2}{(\sum^n_{i=1}\h_i)^{2}} \right] & \ge
& {m_1m_2 \over M^{2}} -2{ m_2\mbox{var}(\h_1)
+m_1\mbox{var}(\h_2)\over M^{3}},
\label{eta12} \\
\bP\left[\frac{\h_1^{2}}{(\sum^n_{i=1}\h_i)^{2}} \right] & \ge &
{\bP[\h_1^{2}] \over M^{2}} \lef( 1+{2m_1 \over M} \ri) -2{ \bP[\h_1^{3}]
\over M^{3}}, \label{eta11} \edmn where $m_i=\bP[\h_i]$ and
$M=\sum_{i=1}^nm_i$.
\end{lemma}
Proof: We set
$$
U=\sum_{i=1}^n(\h_i-m_i)=\sum_{i=1}^n\h_i-M >-M.
$$
Note that $(u+M)^{-2} \ge M^{-2}(1 -{2u \over M})$ for $u \in
(-M,\8)$. Thus, we have that \bdnn \bP\lef[ { \h_1\h_2 \over
(\sum_{i=1}^n\h_i)^{2} } \ri] & = &  \bP\lef[{ \h_1\h_2 \over ( U+M
)^{2}} \ri]
 \ge
M^{-2} \lef( m_1m_2-{2 \over M}\bP\lef[\h_1\h_2U \ri] \ri) \\
\bP\lef[\h_1\h_2U \ri] &= &\bP\lef[\h_1\h_2(\h_1-m_1)
\ri]+\bP\lef[\h_1\h_2(\h_2-m_2) \ri] =m_2 \mbox{var}(\h_1)+m_1
\mbox{var}(\h_2). \ednn These prove (\ref{eta12}). Similarly,
\bdnn \bP\lef[ { \h_12 \over (\sum_{i=1}^n\h_i)^{2} } \ri] & = &
\bP\lef[{ \h_12 \over ( U+M )^{2}} \ri]
 \ge
M^{-2} \lef( \bP\lef[\h_1^{2} \ri]-{2 \over M}\bP\lef[\h_1^2U \ri] \ri),\\
\bP\lef[\h_1^2U \ri] & = &\bP\lef[\h_1^{2}(\h_1-m_1) \ri]
=\bP\lef[\h_1^{3} \ri]-m_1\bP\lef[\h_1^{2} \ri]. \ednn These prove
(\ref{eta11}). \hfill $\Box$

\bigskip
As an immediate consequence, we have (by applying \Lem{CH} to
$\alpha_i \eta_i$ instead of $\eta_i$):
\begin{corollary}
Let $\h_i$, $1 \le i \le n$ ($n \ge 2$) be positive i.i.d.r.v.'s
on a probability space with the probability measure $\bP$, such
that $\bP[\h_i^{3} ] <\8$ for $i=1,..,n$. Then, for any $\alpha_i
\ge0$ satisfying $\sum_{i=1}^n \alpha_i=1$, we have
\begin{eqnarray}
    \bP \Big( { \eta_1 \eta_2 \over \big( \sum_{i=1}^n
    \alpha_i \eta_i\big)^{2} } \Big) & \ge &   1 -
(\bP (\tl{\h}_1^{2})-1)( \alpha_1 +
    \alpha_2) ,   \label{eta13}\\
   \bP \Big( { (\eta_1)^{2}   \over \big( \sum_{i=1}^n \alpha_i
    \eta_i\big)^{2} } \Big)
& \ge & (1+2\a_1 )\bP (\tl{\h}_1^{2}) - 2\a_1\bP (\tl{\h}_1^{3}),
\label{eta14}
\end{eqnarray}
where $\tl{\eta}_1:= \eta_1/m_1$.

\end{corollary}

\bigskip

\begin{lemma}\label{lem10}
Assume $Q(q_{t,x}(0)=0)=1$ and $Q(q_{t,x}(1)=1)<1$. Then,
\bdnl{c_0} \inflim_{t\to\infty} {1\over t} \ln N_t \ge c_0, \qquad
a.s. \edn where $c_0=-\ln Q[\sum_{k \ge 1}k^{-1}q_{t,x}]>0$.
\end{lemma}

{\bf Proof:} For any $(t,y,\nu)$, $K_{t,y}^\nu  $ is independent
of $\cF_t$, hence   $$ P \Big( \big(K_{t, y}^{\nu}\big)^{-1} \,
\big|\, \cF_t\Big) =   P \Big( \big(K_{t, y}^{\nu}\big)^{-1} \,
\Big)=e^{-c_0}.$$
It follows by Jensen's inequality that
\begin{eqnarray*}
  P\left( {1\over N_t} \,\big| \, \cF_{t-1} \right) &=&   P\left( \left[ \sum_y \sum_{\nu=1}^{N_{t-1, y}} K_{t-1,y}^\nu\right]^{-1} \,\big| \, \cF_{t-1} \right)
   \\  &=& {1\over N_{t-1}}\, P\left( \left[ {1\over N_{t-1}}\,\sum_y \sum_{\nu=1}^{N_{t-1, y}} K_{t-1,y}^\nu\right ]^{-1} \,\big| \, \cF_{t-1} \right)
   \\  &\le & {1\over N_{t-1}}\, P\left(  {1\over N_{t-1}}\,\sum_y \sum_{\nu=1}^{N_{t-1, y}}\Big( K_{t-1,y}^\nu\Big)^{-1} \,\big| \, \cF_{t-1} \right)
   \\  &= & {e^{-c_0} \over N_{t-1}}.
   \end{eqnarray*}
Hence $ P\left( {1\over N_t}   \right) \le  e^{-c_0t}$, and
(\ref{c_0}) follows from the Borel-Cantelli lemma.
\hfill$\Box$

\bigskip
We denote by $(\pp_n, n=0, 1,...)$ the semigroup of a simple
symmetric random walk on $\z^d$, namely, $\pp_nf(x):= \sum_y
\pp_n(x,y) f(y)$ where $\pp_n(x,y)$ is the probability that the
random walk starting from $x$ lives at $y$ on the $n$-th step.
Plainly, $\pp_1(x,y)= p(x,y)$. We write $\pp= \pp_1$. Let for any
$z \in \z^d$, $$r_l:= \pp_{2l}(z,z)= \pp_{2l}(0,0)\,\, \sim\, c\,
l^{-d/2}, \qquad l\to\infty.$$ For the sake of notational
convenience, we write $\rh_t(x)\equiv \rh_{t,x}$, so that $\rh_t$
stands for a function on $\zd$.
\begin{lemma}  Assume (\ref{1<K<K}).
For any $(y_1, \nu_1)$ and $(y_2, \nu_2)$, $t\ge1$, we have
\begin{eqnarray}  P \Big( {K_{t, y_1}^{\nu_1} \,  K_{t ,
y_2}^{\nu_2}  \over N_{t+1}^{2}} \,  \big|\, \cF_t\Big) & \ge & { 1
\over N_t^{2}} \, \Big[ (\a -1) 1_{(y_1=y_2)} - c_1 \rh_t(y_1) -c_2
\rh_t(y_2) - {c_2 \over N_t}\Big], \label{kty3}
     \end{eqnarray}
on the event $\{ N_{t,y_1}\wedge N_{t,y_2} \ge 1\}$, where $\a
=\frac{Q[m_{t,x}^{2}]}{m^{2}}$ and $c_1$ and $c_2$ are some constants.
     Consequently,
\begin{eqnarray} P \Big( \rho_{t+1}(y_1) \rho_{t+1}(y_2)  \big|\, \cF_t\Big)
     & \ge &   \big( 1- { c_2\over N_t}) \pp \rho_t(y_1) \, \pp \rho_t(y_2)
+ (\a -1) \sum_{z}   \rho_t^{2}(z) p(z ,
     y_1) p(z , y_2)   \nonumber\\
      && \quad - c_1 \Big[ \pp \rho_t(y_1) \, \pp (\rho_t^{2})(y_2) + \pp \rho_t(y_2) \, \pp (\rho_t^{2})(y_1) \Big]
       \nonumber\\
       && \quad- {\a \over N_t} \sum_z p(z, y_1) p(z, y_2)  \rho _t(z). \label{rhoy12}
     \end{eqnarray}
\end{lemma}
{\bf Proof:}
  Firstly, we consider (\ref{kty3}) in the case
$(y_1,\nu_1) \not = (y_2, \nu_2)$. Let $A \in \cF_t$ and $A \sub
\{ N_{t,y_1}\wedge N_{t,y_2} \ge 1\}$. Under $P^q$, $\{K_{t,
y}^\nu\}_{t, \nu}$ are independent (but not identically
distributed) and independent of $1_A, N_{t,\cdot}$. Write
$M_t=\sum_y m_{t,y} N_{t,y}$. Noting that $N_{t+1}=\sum_y \,
\sum_{\nu=1}^{N_{t,y}} K_{t, y}^\nu $ and applying (\ref{eta12})
to $\eta_1=K_{t, y_1}^{\nu_1}$ and $\eta_2=K_{t , y_2}^{\nu_2}$,
we get
$$
P^q \Big( 1_A{ K_{t, y_1}^{\nu_1} \,  K_{t , y_2}^{\nu_2}  \over
N_{t+1}^{2}} \,\Big)
\ge P^q \lef(  1_A  { m_{t,y_1}m_{t,y_2} \over M_t^{2}}\rig) - 2P^q
\lef(1_A{m_{t,y_2}  m^{(2)}_{t,y_1} +m_{t,y_1}  m^{(2)}_{t,y_2}
\over N_t^{3}}\rig),
$$
since $M_t \ge N_t$. Therefore, by taking $Q$-expectation,
$$
P \Big( 1_A{ K_{t, y_1}^{\nu_1} \,  K_{t , y_2}^{\nu_2}  \over
 N_{t+1}^{2}} \,\Big)
\ge P \lef(  1_A  { m_{t,y_1}m_{t,y_2} \over M_t^{2}}\rig) - 2P
\lef(1_A{m_{t,y_2}  m^{(2)}_{t,y_1} +m_{t,y_1}  m^{(2)}_{t,y_2}
\over N_t^{3}}\rig).
$$
Observe that under $\P$, $m_{t,\cdot}$ are i.i.d. and independent
of $\cF_t$. It turns out from (\ref{eta13}) and (\ref{eta14}) that
\bdnn P \lef( { m_{t,y_1}m_{t,y_2} \over M_t^{2}} \big|\, \cF_t\rig)
& = & {1 \over N_t^{2}} P\lef( { m_{t,y_1}m_{t,y_2} \over (
 \sum_y \rho_t(y)  m_{t,y} )^{2}}\, \big|\, \cF_t\right) \\
& \ge & {1 \over N_t^{2}}\lef[ 1 + (\a -1)
 1_{(y_1=y_2)} - c_1 \rh_t(y_1) -c_2 \rh_t(y_2) \rig].
\ednn On the other hand, we have
$$
P\lef(  m_{t,y_2}  m^{(2)}_{t,y_1} +m_{t,y_1}  m^{(2)}_{t,y_2}|\,
\cF_t \rig) \le 2m^{(3)}<\8
$$
by our integrability assumption. Hence, with $c_2=4m^{(3)}$,
$$ P \Big( 1_A \, {K_{t, y_1}^{\nu_1}
\,  K_{t , y_2}^{\nu_2}  \over N_{t+1}^{2}} \, \Big)  \ge P\left(
{1_A \over N_t^{2}} \, \Big[ 1+ (\a -1) 1_{(y_1=y_2)} - c_1
\rh_t(y_1) -c_2 \rh_t(y_2) - {c_2 \over N_t}\Big]\right),
$$ yielding
 (\ref{kty3}) in the case $(y_1, \nu_1) \not= (y_2, \nu_2)$. The
 case $(y_1, \nu_1)   = (y_2, \nu_2)$ is obtained in the same way by
 applying (\ref{eta11}) instead of (\ref{eta12})  and by eventually
 modifying the constants.

 To obtain (\ref{rhoy12}), we have that \begin{eqnarray} P \Big( \rho_{t+1}(y_1) \rho_{t+1}(y_2)  \big|\, \cF_t\Big)
     & =& \sum_{z_1, z_2} \sum_{\nu_1=1}^{N_{t, z_1}}
     \sum_{\nu_2=1}^{N_{t,z_2}} P\left( { \delta_{y_1}( X^{\nu_1}_{t, z_1})
\delta_{y_2}( X^{\nu_2}_{t,     z_2}) K_{t, z_1}^{\nu_1} K_{t,
z_2}^{\nu_2} \over  N_{t+1}^{2}} \,
\big|\, \cF_t \right) \nonumber \\
     &=& \sum_{z_1, z_2} \sum_{\nu_1=1}^{N_{t,z_1}}
     \sum_{\nu_2=1}^{N_{t,z_2}} h_{1,2}\,  P\left( { K_{t, z_1}^{\nu_1} K_{t, z_2}^{\nu_2} \over  N_{t+1}^{2}} \, \big|\, \cF_t \right)
     \label{h11}\end{eqnarray}

     by means of the
     independence between $(X^{\nu_1}_{t, z_1}, X^{\nu_2}_{t,
     z_2})$ and $(\cF_t, N_{t+1},K_{t, z_1}^{\nu_1} K_{t, z_2}^{\nu_2})$,
     and the function $h_{1,2}$ is defined as follows: \begin{eqnarray} h_{1,2}&:= &P \Big(\delta_{y_1}( X^{\nu_1}_{t, z_1}) \delta_{y_2}( X^{\nu_2}_{t,
     z_2}) \Big)  \nonumber \\
     &= &1_{((z_1, \nu_1)=(z_2, \nu_2))} p(z_1, y_1)
     1_{(y_1=y_2)} + 1_{((z_1, \nu_1)\not=(z_2, \nu_2))} p(z_1,
     y_1) p(z_2, y_2)  \label{h12}\\
     &\ge &  1_{((z_1, \nu_1)\not=(z_2, \nu_2))} p(z_1,
     y_1) p(z_2, y_2). \nonumber\end{eqnarray}

Applying (\ref{kty3}) we get \begin{eqnarray*} P \Big(
\rho_{t+1}(y_1) \rho_{t+1}(y_2)  \big|\, \cF_t\Big)
     & \ge & \sum_{z_1, z_2} \sum_{\nu_1, \nu_2} h_{1,2} \,
     {1\over N_t^{2}} \, \Big[ 1+ (\a -1) 1_{(z_1=z_2)}
- c_1 \rh_t(z_1) -c_2 \rh_t(z_2) - {c_2 \over N_t}\Big] \\
  &\ge & \sum_{(z_1, \nu_1)\not=(z_2, \nu_2)} p(z_1,
     y_1) p(z_2, y_2) \,
     {1\over N_t^{2}} \, \Big[ g_t(z_1, z_2) + (\a -1) 1_{(z_1=z_2)}
     \Big],
\end{eqnarray*} with $g_t(z_1, z_2)= 1  - c_1 \rh_t(z_1) -c_2 \rh_t(z_2)
- {c_2 \over N_t}.$ Let us compute explicitly the
 above sum $\sum_{(z_1, \nu_1)\not=(z_2, \nu_2)}\cdot\cdot\cdot$:
\begin{eqnarray*}
 \sum_{(z_1, \nu_1)\not=(z_2, \nu_2)}&=& \sum_{z_1\not= z_2}
 N_{t,z_1} N_{t,z_2}  p(z_1,
     y_1) p(z_2, y_2) \,
     {1\over N_t^{2}} \,  g_t(z_1, z_2)   \\
     &\qquad & + \sum_{z} ( N_{t,z}^{2}-
     N_{t,z})  p(z ,
     y_1) p(z , y_2) \,
     {1\over N_t^{2}} \,  \big[ g_t(z, z) +\a -1\big],
 \end{eqnarray*}

\noindent by removing  the diagonal terms. Using the definition of
$\rho_t(z)= N_{t,z}/N_t$, we get
\begin{eqnarray*}
 \sum_{(z_1, \nu_1)\not=(z_2, \nu_2)}&=& \sum_{z_1, z_2}
 \rho_t(z_1) \rho_t(z_2)  p(z_1,
     y_1) p(z_2, y_2)  \,  g_t(z_1, z_2)  \\
     \quad && + (\a -1) \sum_{z}   \rho_t2(z) p(z ,
     y_1) p(z , y_2)   \, - \sum_{z}
      p(z ,
     y_1) p(z , y_2) \,
     {\rho_t(z)\over N_t} \,  \big[ g_t(z, z) +\a -1\big] \\
     &\ge & \big( 1- { c_2\over N_t}) \pp \rho_t(y_1) \, \pp \rho_t(y_2)
+ (\a -1) \sum_{z}   \rho_t2(z) p(z ,
     y_1) p(z , y_2)  \\
      && \,\, - c_1 \Big[ \pp \rho_t(y_1) \, \pp (\rho_t^{2})(y_2) + \pp \rho_t(y_2) \, \pp (\rho_t^{2})(y_1) \Big]
       - {\a \over N_t} \sum_z p(z, y_1) p(z, y_2)  \rho _t(z),
     \end{eqnarray*}
as desired.
 \hfill$\Box$.

\bigskip

Recall that $\R_t=\sum_x \rho_t^{2}(x) $.  Let $t \ge 2$. The
following lemma shows the r\^ole de semigroup in $\R_t$.

\begin{lemma} \label{lem20}
Assume (\ref{1<K<K}). There exists a constant $c_3>0$ such that
for all $1\le s \le t-1$,
\begin{eqnarray*}
  P \left( \sum_x \big( \pp_{t- (s+1)} \, \rho_{s+1}(x)\big)^{2} \,
  \big|\, \cF_s \right) & \ge & \sum_x \big(
\pp_{t-s} \rho_s(x)\big)^{2} + (\a -1) \, r_{t-s} \R_s - 2 c_1
\R_s^{3/2} - {c_3\over N_s}.
\end{eqnarray*}
\end{lemma}

{\bf Proof:} Observe that $$ \sum_x \big( \pp_{t- (s+1)} \,
\rho_{s+1}(x)\big)^{2} = \sum_x \sum_{y_1, y_2} \pp_{t- (s+1)}(x,
y_1) \pp_{t- (s+1)}(x, y_2) \rho_{s+1}(y_1) \rho_{s+1}(y_2).$$

Applying (\ref{rhoy12}) gives $$  P \left( \sum_x \big( \pp_{t-
(s+1)} \, \rho_{s+1}(x) \big)^{2} \,
  \big|\, \cF_s \right) \ge \big( 1- { c_2\over N_s}) I_2 + (\a -1)
  I_3 - c_1 I_4 - {\a \over N_s} I_5, $$

  \noindent with \begin{eqnarray*}
    I_2 &:=& \sum_x \sum_{y_1, y_2} \pp_{t- (s+1)}(x,
y_1) \pp_{t- (s+1)}(x, y_2) \pp \rho_s(y_1) \, \pp \rho_s(y_2),
   \\ I_3&:=& \sum_x \sum_{y_1, y_2} \pp_{t- (s+1)}(x,
y_1) \pp_{t- (s+1)}(x, y_2) \sum_z \rho_s^{2}(z) p(z ,
     y_1) p(z , y_2) ,
    \\ I_4 &:=& \sum_x \sum_{y_1, y_2} \pp_{t- (s+1)}(x,
y_1) \pp_{t- (s+1)}(x, y_2) \Big[ \pp \rho_s(y_1) \, \pp
(\rho_s^{2})(y_2) + \pp \rho_s(y_2) \, \pp (\rho_s^{2})(y_1) \Big],
    \\ I_5&:=& \sum_x \sum_{y_1, y_2} \pp_{t- (s+1)}(x,
y_1) \pp_{t- (s+1)}(x, y_2) \sum_z p(z, y_1) p(z, y_2)  \rho
_s(z).
  \end{eqnarray*}

Using the semigroup property and noting that
$\sum_x\big(\pp_{t-s}(x, z) \big)^{2} = \pp_{2t-2s}(z, z) = r_{t-s} $
for any $z$, we obtain
\begin{eqnarray*} I_2 &=& \sum_x \big(
\pp_{t-s} \rho_s(x)\big)^{2}, \\
    I_3&=& \sum_x \sum_z \big( \pp_{t-s}(x, z)\big)^{2} \rho^{2}_s(z)
    = \sum_z \sum_x\big( \pp_{t-s}(x, z)\big)^{2} \rho^{2}_s(z) =
    r_{t-s} \, \sum_z \rho^{2}_s(z) , \\
    I_4&=& 2 \sum_x \, \pp_{t-s} \rho _s(x) \, \pp_{t-s}(\rho
    ^2_s) (x), \\
    I_5&=& \sum_x \sum_z \big( \pp_{t-s}(x,z)\big)^{2} \rho_s(z)=
    \sum_z \rho _s(z) r_{t-s}= r_{t-s}.
\end{eqnarray*}
By the translation invariance and Cauchy-Schwarz's inequality, we
see that
$$
\R_s= \sum_x \pp_{t-s} (\rho ^2_s)(x) \ge \sum_x \big(\pp_{t-s}
(\rho _s)(x)\big)^{2} \ge \max_x \pp_{t-s} (\rho _s)(x)^{2},
$$
and hence that $I_4 \le 2 \, \R_s^{3/2}$. This implies the lemma
with $c_3=\a + c_2$.
\hfill $\Box$

\vvs Define $$ V_t= \sum_{s=1}^t \, \R_s, \qquad t=1,2,...$$
\begin{lemma} \label{lem30}
  Assume (\ref{1<K<K}).
Fix  $j \ge 0$. The martingale $Z_j(\cdot)$ defined by $$ Z_j(t):
=\sum_{s=1}^t \left( \sum_x \big( \pp_j \rho _s(x)\big)^{2} - P \Big(
\sum_x \big( \pp_j \rho _s(x)\big)^{2} \, \big|\, \cF_{s-1}\Big)
\right), \qquad t \ge 1.$$ satisfies the following law of large
numbers:
$$
\{ V_\8 =\8\} \st{\mbox{\scriptsize a.s.}}{\sub} \{ { Z_j(t) \over
V_t} \, \to \, 0, \qquad t \to\infty, \, \}.$$
\end{lemma}
{\bf Proof:} Let us compute the increasing process $\langle
Z_j\rangle_\cdot $ associated to $Z_j$. By Cauchy-Schwarz'
inequality, $ \big( \sum_x \pp_j \rho _s(x)\big)^{2} \le  \sum_x
\pp_j  \rho^{2} _s(x) = \R_s\le 1$. It follows that
\begin{eqnarray*}
\big( Z_j(s) - Z_j(s-1)\big)^{2} &\le & 2 \left( \sum_x \big( \pp_j
\rho _s(x)\big)^{2}  \right)^{2} + 2 \left( P \Big( \sum_x \big( \pp_j
\rho _s(x)\big)^{2} \, \big|\,
\cF_{s-1}\Big) \right)^{2} \\
    &\le & 2 \R^{2}_s + 2 P \left( \R_s \,\big|\,
\cF_{s-1}\ri)^{2} \\
    &\le & 2 \R_s + 2 P\Big( \R_s \,\big|\,
\cF_{s-1}\Big).
\end{eqnarray*}
Hence,
$$ \langle Z_j\rangle_s - \langle Z_j\rangle_{s-1}=
P\left( \big( Z_j(s) - Z_j(s-1)\big)^{2} \,\big|\, \cF_{s-1} \right)
\le 4 P\Big( \R_s \,\big|\, \cF_{s-1}\Big).$$
We will prove that
\begin{equation}\label{c7}
P\Big( \R_s \,\big|\, \cF_{s-1}\Big) \le 2m^{(2)} \R_{s-1}.
\end{equation}
Then,
 $\lan Z_j\ran_t \le 8m^{(2)} V_{t-1}$, and the
lemma follows from  the standard law of large numbers for a
square-integrable martingale, c.f. section \ref{pf_loc},(4).

It remains to show (\ref{c7}).   Using (\ref{h11}) and (\ref{h12})
to $y_1=y_2=y$, we have
\begin{eqnarray*} P\Big( \R_s \,\big|\, \cF_{s-1}\Big) &=& \sum_y
 P \Big( \rho^{2} _s(y) \,\big|\, \cF_{s-1}\Big) \\
    &=& \sum_y \sum_{z_1, z_2} \sum_{\nu_1=1}^{N_{s-1,z_1}}
     \sum_{\nu_2=1}^{N_{s-1,z_2}} h_{1,2}\,
P\left( { K_{s-1, z_1}^{\nu_1} K_{s-1, z_2}^{\nu_2} \over  N_s^{2}}
\, \big|\, \cF_{s-1} \right)
    \\ &\le & \sum_y \sum_{z_1, z_2} \sum_{\nu_1=1}^{N_{s-1,z_1}}
     \sum_{\nu_2=1}^{N_{s-1,z_2}} h_{1,2}\,  {m^{(2)}\over N^{2}_{s-1}}.
    \end{eqnarray*}
To obtain the last inequality, we used $N_s \ge N_{s-1}$ and the
independence between $K_{s-1, \cdot}^{\cdot}$ and $\cF_{s-1}$.
We divide the last summation into bound the summation over
$(z_1,\n_1)=(z_2,\n_2)$ and that over $(z_1,\n_1)=(z_2,\n_2)$, to
see that
$$
\sum_y \sum_{z_1, z_2} \sum_{\nu_1=1}^{N_{s-1,z_1}}
     \sum_{\nu_2=1}^{N_{s-1,z_2}}
{h_{1,2}\over N^{2}_{s-1}} \le { 1 \over N_{s-1}} +
     \sum_x \big( \pp \rho _{s-1}(x)\big)2 \le { 1 \over N_{s-1}}
     + \R_{s-1} \le 2\R_{s-1}.$$
Here, we used $\R_{s-1}= \sum_x N_{s-1,x}^{2} /N^{2}_{s-1} \ge
1/N_{s-1}$ to see the last inequality. Putting things together, we
 have  (\ref{c7}) and the proof of the lemma is now complete.
\hfill$\Box$
\subsection{Proof of \Thm{sloc}:}
We first note that there are $\epsilon >0$ and $t_0 \in \N$ such
that \bdnl{eT_0} \sum_{s=1}^{t_0} r_s \ge {1+\epsilon \over \a
-1}. \edn For $d=1,2$, we take $\epsilon =1$. Then, (\ref{eT_0})
holds for $t_0$ large enough, since $\sum_{s=1}^{\8} r_s=\8$. For
$d \ge 3$, our assumption  $P(\ov{N}_\8 =0)=1$ implies $\a \ge
\a^* >1/\pi_d$ by \Prop{WD}. Since $\sum_{s=1}^\8
r_s=\frac{\pi_d}{1-\pi_d}$, as is well known, (\ref{eT_0}) holds
for small enough $\epsilon >0$ and large enough $t_0$.

Let $t> t_0$. Applying Lemma \ref{lem20} to $s=t-1, t-2,..., t-
t_0$ and taking the sum on $s$, we get  \begin{eqnarray*}
   & & \sum_{s=t-t_0}^{t-1} (2 c_1   \R_s^{3/2} + {c_3\over N_s}) \\
   &\ge&
  \sum_{s=t-t_0}^{t-1} \left(  \sum_x \big(
\pp_{t-s} \rho_s(x)\big)^{2}- P \Big( \sum_x \big( \pp_{t- (s+1)} \,
\rho_{s+1}(x)\big)^{2} \,  \big|\, \cF_s \Big)  \right)
 + (\a -1)\,\sum_{s=t-t_0}^{t-1}r_{t-s} \R_s \\
    &=& \sum_{s=t-t_0}^{t-1}
\left( \sum_x \big( \pp_{t- (s+1)} \, \rho_{s+1}(x)\big)^{2} - P
\Big( \sum_x \big( \pp_{t- (s+1)} \, \rho_{s+1}(x)\big)^{2} \,
  \big|\, \cF_s \Big)   \right) +
      \\ && \quad\sum_{s=t-t_0}^{t-1} \left(   \sum_x \big(
\pp_{t-s} \rho_s(x)\big)^{2} - \sum_x \big( \pp_{t- (s+1)} \,
\rho_{s+1}(x)\big)^{2}\right)  + (\a -1) \,\sum_{s=t-t_0}^{t-1}
r_{t-s} \R_s \\
    &=& \sum_{s=t-t_0}^{t-1} \Big[ Z_{t-(s+1)}(s+1) -
    Z_{t-(s+1)}(s)\Big] + \sum_x \big( \pp_{t_0} \,
\rho_{t-t_0}(x)\big)^{2} - \R_t + (\a -1) \,\sum_{s=t-t_0}^{t-1}
r_{t-s} \R_s   ,
\end{eqnarray*}

where we recall that the   martingale $Z_j(\cdot)$ are defined in
Lemma \ref{lem30}. By change of variable $s=t-j$, we have proven
that
$$ \sum_{j=1}^{t_0} (2 c_1   \R_{t-j}^{3/2} + {c_3\over N_{t-j}})
\ge \sum_{j=1}^{t_0} \Big[ Z_{j-1}(t-j+1) - Z_{j-1}(t-j)\Big]   -
\R_t + (\a -1)\,\sum_{j=1}^{t_0} r_j \R_{t-j}.$$

Taking the sum of these inequalities for $t=t_0+1,..., T$, we
obtain that \begin{eqnarray*}  \sum_{t=t_0+1}^T\sum_{j=1}^{t_0} (2
c_1 \R_{t-j}^{3/2} + {c_3\over N_{t-j}}) &\ge& \sum_{j=1}^{t_0}
\Big[ Z_{j-1}(T-j+1) - Z_{j-1}(t_0-j+1)\Big] - ( V_T - V_{t_0})
\\ && \qquad + (\a -1) \,\sum_{j=1}^{t_0} r_j \big( V_{T-j} - V_{t_0-j}\big).
\end{eqnarray*}

Since $\R_s \le 1$, \bdnn V_{T-j} - V_{t_0-j} & \ge & V_T -j -
(t_0-j)= V_T
-t_0, \\
 (\a -1)\sum_{j=1}^{t_0} r_j \big( V_{T-j} -V_{t_0-j}\big)
& \ge & (\a -1)\sum_{j=1}^{t_0} r_j  V_T - c_9 \ge (1+\epsilon
)V_T - c_9, \ednn
 with constant $c_9= (\a -1) t_0 \sum_{j=1}^{t_0} r_j.$
Hence,
\begin{equation}\label{eq50} \sum_{t=t_0+1}^T\sum_{j=1}^{t_0} (2 c_1
\R_{t-j}^{3/2} + {c_3\over N_{t-j}}) \ge \sum_{j=1}^{t_0} \Big[
Z_{j-1}(T-j+1) - Z_{j-1}(t_0-j+1)\Big]  + \epsilon V_T - c_9.
\end{equation}

Recall from Lemma \ref{lem10} that $ \sum_{t=1}^\infty {1 \over
N_t}\, < \infty, a.s.,$ which combined with   Lemma   \ref{lem30}
imply that the two sums involving respectively  ${c_3\over
N_{t-j}}$ and $Z_{j-1}(T-j+1)$ in (\ref{eq50}) are negligible,
 relative to  $V_T$. It follows that
$$ \liminf_{T\to\infty} {1 \over V_T}
\sum_{t=t_0+1}^T\sum_{j=1}^{t_0} \R_{t-j}^{3/2}
 \ge  { \epsilon\over 2
c_1} , \quad a.s.$$

Consequently, $$ \liminf_{T\to\infty} {1 \over V_T} \sum_{t= 1}^T\
\R_t^{3/2}   \ge  { \epsilon\over 2 c_1 \, t_0} , \quad a.s.,$$

which implies that $$ \limsup_{t \to\infty} \R_t \ge \big({
\epsilon \over 2 c_1 \, t_0}\big)^{2}, \qquad a.s.$$  This completes
the proof of Theorem. \hfill $\Box$

%



{\footnotesize

\baselineskip=12pt

\noindent
\begin{tabular}{lll}
&\hskip20pt Yueyun Hu
    & \hskip45pt Nobuo Yoshida\\
&\hskip20pt D\'epartement de Math\'ematiques
    & \hskip45pt Division of Mathematics, Graduate School of Science\\
&\hskip20pt Universit\'e Paris XIII
    & \hskip45pt Kyoto University\\
&\hskip20pt F-93430 Villetaneuse
    & \hskip45pt Kyoto 606-8502\\
&\hskip20pt France
    & \hskip45pt Japan \\
&\hskip20pt {\tt yueyun@math.univ-paris13.fr}
    & \hskip45pt {\tt nobuo@math.kyoto-u.ac.jp}
\end{tabular}

\end{document}